\documentclass[11pt]{amsart}
\usepackage{url}
\usepackage{inputenc}
\usepackage{latexsym,color,amsmath,amsthm,amssymb,amscd,amsfonts,bbm}
\usepackage{graphicx}
\usepackage{scalefnt}
\usepackage[font=small,labelfont=bf]{caption}
\usepackage{amsopn}
\usepackage{relsize}
\usepackage{hyperref}
\usepackage{mathrsfs}
\usepackage[dvipsnames]{xcolor}

\newtheoremstyle{nonum}{}{}{\itshape}{}{\bfseries}{.}{ }{\thmnote{#3}}

\newtheorem{thm}{Theorem}
\newtheorem*{thm*}{Theorem}

\newtheorem*{cor*}{Corollary}
\newtheorem{lem}[thm]{Lemma}
\newtheorem*{lem*}{Lemma}

\newtheorem*{rem*}{Remark}
\newtheorem{conj}[thm]{Conjecture}
\newtheorem*{conj*}{Conjecture}

\newtheorem*{definition*}{Definition}

\newtheorem*{fact*}{Fact}

\newtheorem*{rems*}{Remarks}

\newtheorem{prob}[thm]{Problem}
\newtheorem*{prob*}{Problem}

\theoremstyle{nonum}

\newcommand{\R}{\mathbb R}
\newcommand{\RR}{\mathbb R}

\newcommand{\Sph}{\mathbb S}

\newcommand{\rn}{\mathbb R^n}
\newcommand{\sn}{\mathbb S^{n-1}}
\newcommand{\chara}[1]{{\mathbbm{1}_{#1}}}

\def\Vol{{\rm Vol}}

\newcommand{\iprod}[2]{\langle #1,#2 \rangle} 

\newcommand{\kpn}{\mathcal{K}_+^n}
\newcommand{\Hn}{\mathcal{H}^{n-1}}

\newcommand{\wt}{\widetilde}

\newcommand{\barV}{\bar{V}}

\makeatletter
\def\moverlay{\mathpalette\mov@rlay}
\def\mov@rlay#1#2{\leavevmode\vtop{%
		\baselineskip\z@skip \lineskiplimit-\maxdimen
		\ialign{\hfil$\m@th#1##$\hfil\cr#2\crcr}}}
\newcommand{\charfusion}[3][\mathord]{
	#1{\ifx#1\mathop\vphantom{#2}\fi
		\mathpalette\mov@rlay{#2\cr#3}
	}
	\ifx#1\mathop\expandafter\displaylimits\fi}
\makeatother

\begin{document}
	\title{Brunn-Minkowski and Reverse Isoperimetric Inequalities for Dual Quermassintegrals}
	\author{Shay Sadovsky}
	\address{Courant Institute of Mathematical Sciences, New York University, NY}
	\email{ss20011@nyu.edu}
	
	\author{Gaoyong Zhang}
	\address{Courant Institute of Mathematical Sciences, New York University, NY}
	\email{gaoyong.zhang@nyu.edu}

\subjclass{52A40}

\keywords{Convex body, cube, dual Brunn-Minkowski theory, dual quermassintegrals, 
Brunn-Minkowski inequality, John's position, reverse isoperimetric inequality.}

\thanks{Research of Zhang supported, in part, by NSF Grant  DMS--2005875.}

	\maketitle

\begin{abstract}
	This paper establishes two new geometric inequalities in the dual Brunn-Minkowski theory. The first, originally conjectured by Lutwak, is the Brunn-Minkowski inequality for dual quermassintegrals of
origin-symmetric convex bodies. The second, generalizing Ball's volume ratio inequality, is a reverse isoperimetric inequality: among all origin-symmetric convex bodies in John's position, the cube maximizes the dual quermassintegrals.
\end{abstract}

\section{Introduction}
In his seminal paper \cite{lutwak1975}, Erwin Lutwak introduced the dual mixed volumes as a natural dual notion to that of the classical mixed volume (due to Minkowski, for more see \cite{Schneider-book-NEW}). These dual mixed volumes are defined there as follows: Given $K_1,\dots K_n \subseteq \RR^n$ compact convex sets with non-empty interiors that contain the origin, henceforth referred to as convex bodies, their dual mixed volume is given by
\begin{equation}
	\wt{V} (K_1,\dots, K_n)=\frac{1}{n} \int_{\Sph^{n-1}} \rho_{K_1}(u) \dots \rho_{K_n}(u)\, du
\end{equation}
where $\Sph^{n-1}$ is the usual $n-1$ dimensional unit sphere, and $\rho_K(x) = \sup\{\lambda \ge 0: \lambda x \in K\}$ is the radial function of the convex body $K$.

In the case where $K_1,\dots, K_j = K$ are equal, and $K_{j+1},\dots,K_n = B_2^n$ are equal to the unit ball, we denote
\begin{equation}\label{eq:def-radial}
	\wt{V}_{j}(K) = \wt{V} (K_1,\dots, K_n) = \frac{1}{n} \int_{\Sph^{n-1}}  \rho_{K}^j(u)\, du,
\end{equation}
and call this quantity the $j$-th dual quermassintegral or $j$-th dual volume of $K$. Indeed, this name is fitting, as the dual quermassintegrals have a geometric formulation as integrals over volumes of lower dimensional sections,
$$\wt{V_j}(K)=\frac{\omega_n}{\omega_j}\int_{G_{n,j}} \Vol_j(K\cap \xi)d\xi$$
in analogy to the definition of the usual quermassintegrals as integrals over volume of lower dimensional projections, where $\omega_j$ denotes the volume of the $j$-dimensional unit ball
and $G_{n,j}$ is the Grassmann manifold of $j$-dimensional subspaces.

For any $q\in \R$, the $q$-th dual quermassintegral of the convex body $K$ is defined by
\begin{equation}\label{d-v}
\wt V_q(K) = \frac1n \int_{\sn} \rho_K^q (u)\, du.
\end{equation}

We state the following well known alternative definition of the dual quermassintegral (as can be seen e.g. in \cite{gardner2007dual}). For $K\subseteq \RR^n$
\begin{equation}\label{eq:def-norm}
	\wt{V}_{q}(K) = \frac qn \int_K |x|^{q-n} \, dx, \ \ \ q>0.
\end{equation}
One may check that \eqref{eq:def-radial} and \eqref{eq:def-norm} are equivalent by switching to polar coordinates.

The normalized $q$-th dual quermassintegral of $K$ is
\begin{equation}\label{ndv}
\barV_q(K) = \Big(\frac1{n\omega_n} \int_{\sn} \rho_K^q (u)\, du\Big)^\frac1q = \Big(\frac{\wt V_q(K)}{\omega_n}\Big)^\frac1q.
\end{equation}
Then
\[
\barV_0(K) = e^{\frac1{n\omega_n}E(K)}, \ \ \ E(K)=\int_{\sn} \log \rho_K(u)\, du,
\]
where $E(K)$ is called the dual entropy of $K$.

In \cite{lutwak1975}, Lutwak proved the elementary properties of the dual mixed volumes, as well as several inequalities analogous to fundamental inequalities in convex geometry, such as a dual Alexandrov-Fenchel inequality and a dual isoperimetric inequality. These examples are only a fraction of the many dual problems and properties which can be found in what is today called the ``dual Brunn-Minkowski theory". In particular there exists a radial sum analogous to  the Minkowski sum, a dual Alexandrov-Fenchel inequality, a dual isoperimetric inequality, and a dual Minkowski problem, introduced in \cite{HLYZ-2018} analogously to the usual Minkowski problem of geometric measures. We refer the reader to the extensive work of Huang, Lutwak, Yang and Zhang on the dual Brunn-Minkowski theory in \cite{HLYZ-2016geometric} and the recent survey \cite{HYZ-2025-survey}.

In this work, we extend the dual theory to two new inequalities. The first of these inequalities is the Brunn-Minkowski inequality  of dual quermassintegrals which has been conjectured by Lutwak and presented in various talks and seminars by the group Lutwak-Yang-Zhang and their collaborators (as stated, e.g. in \cite{HYZ-2025-survey}).

\begin{conj}[Lutwak]
	If $K$ and $L$ are origin-symmetric convex bodies, then
	\[
	\wt V_j (K+L)^\frac1j \ge \wt V_j (K)^\frac1j + \wt V_j (L)^\frac1j, \ \ \ j=2, \ldots, n-1,
	\]
	with equality if and only if $K$ and $L$ are dilates of each other.
\end{conj}

We begin with showing that the conjectured Brunn-Minkowski inequality for $\wt{V}_j$ is true.
Recently, an amazing breakthrough in proving general Brunn-Minkowski inequalities by
using analytic tools has been achieved by
 Kolesnikov-Milman \cite{KM-ajm2018, KM-jfa2017} and
  Kolesnikov-Livshyts \cite{kolesnikov-livshyts}, which allowed Cordero-Erausquin and Rotem in \cite{cordero-rotem} to prove  Brunn-Minkowski inequalities for a large class of measures in $\rn$.
  As a corollary of these recent works, we present a solution to the conjecture, extended to all $q> 0$.

\begin{thm}\label{thm:BM}
	Let $K,L$ be origin-symmetric convex bodies in $\RR^n$. Then, for $0< q \le n$,
	\begin{equation}\label{eq:BM-Wj}
		 \wt{V}_{q}(K+ L)^{\frac{1}{q}} \ge \wt{V}_{q}(K)^{\frac{1}{q}}
 + \wt{V}_{q}(L)^{\frac{1}{q}},
	\end{equation}
\end{thm}
with equality if $K$ and $L$ are dilates of each other.

The equality cases of Lutwak's conjecture have yet to be resolved, and in fact are only known for the case of $q\le 1$, for which the inequality together with its equality condition were shown by Xi-Zhang in \cite{xi-zhang}. Note that the case of $q=n$ of \eqref{eq:BM-Wj}
is the classical Brunn-Minkowski inequality of origin-symmetric convex bodies.

We present the results of \cite{kolesnikov-livshyts} and \cite{cordero-rotem} and the resulting proof of the theorem in Section \ref{sec:dual-BM}.

Our second theorem follows the path of the reverse isoperimetric inequalities, firstly due to Ball \cite{Ball-lnm1988, ball1991volume} and later to Barthe \cite{barthe1998extremal}, Schechtman-Schmuckenschl\"ager \cite{schechtman1995concentration}, Schmuckenschl\"ager \cite{schmuckenschlager1999extremal},
Lutwak-Yang-Zhang \cite{LYZ-JDG2004, LYZ-AJM2007, LYZ-JDG2010}, Li-Leng \cite{LL-MZ2012},
Schuster-Weberndorfer \cite{SW-jdg2012}, Li-Huang-Xi \cite{LHX-aam2017}, and others.
They showed that several quantities of convex bodies are extremized by cubes, cross-polytopes,
 or simplices when the bodies are restricted to either John or Lowner position. We present
 such a reverse isoperimetric inequality for dual
quermassintegrals.

Denote by $B_p^n = \{x\in \RR^n:\|x\|_p\le 1\}$, $1\le p\le \infty$.

\begin{thm}\label{thm:reverse-isop}
	Let $K\subseteq \RR^n$ be an origin-symmetric convex body.
	If $K$ is in John's position, i.e. if $B_2^n$ is the ellipsoid of largest volume  contained in $K$, then
	\begin{equation}\label{eq:thm-reverse-isop}
		\barV_q(K) \le \barV_q(B_\infty^n), \ \  q\in (-\infty, n],    
      \end{equation}
      with equality in each of the inequalities if and only if $K$ is a rotation of $B_\infty^n$.
\end{thm}

The case of $q=n$ of \eqref{eq:thm-reverse-isop} is Ball's volume ratio inequality of origin-symmetric
convex bodies \cite{Ball-lnm1988}.  For an arbitrary convex body $K$ that is in John's position,
Ball's volume ratio inequality says that the volume of $K$ is less than the volume of a regular simplex,
see \cite{ball1991volume}.


\section{Brunn-Minkowski inequality for dual quermassintegrals}\label{sec:dual-BM}

In this section, we follow in the footsteps of \cite{KM-jfa2017,  CLM-jfa2017, KM-ajm2018, kolesnikov-livshyts, cordero-rotem} who have proved several extensions of the Brunn-Minkowski inequality using Pointcar\'e type inequalities on weighted manifolds.  We begin with some notation.

Let $\mu$ be a Borel measure in $\rn$. Measure $\mu$ is called log-concave on origin-symmetric convex bodies
if
\begin{equation}\label{l-c}
	\mu((1-\lambda) K + \lambda L) \ge \mu(K)^{1-\lambda} \mu(L)^{\lambda}, \ \  0\le \lambda \le 1,
\end{equation}
for any origin-symmetric convex bodies $K$ and $L$.
Let $\kpn$ be the class of origin-symmetric convex bodies with $C^2$ boundary of positive curvature.
Let $C_e^2(\sn)$ be the set even $C^2$ functions on $\sn$.
For any $\psi \in C_e^2(\sn)$ and $K\in \kpn$, the convex body $K_s$ with support function
\[
h_{K_s} = h_K + s\psi
\]
is also in $\kpn$ when $s$ is small. Thus,
\[
C_e^2(\sn) = \{h_L - h_K : K, L \in \kpn\}.
\]

Denote by $\nu$ the Gauss map of $K\in\kpn$, and $\varphi=\psi(\nu)$. Let $\Hn$ be $(n-1)$-dimensional
Hausdorff measure. Its restriction to the boundary $\partial K$ is denoted by $\Hn_{\partial K}$.

The following lemmas are similar to Lemmas 3.1 and 3.2 in Colesanti-Livshyts-Marsiglietti \cite{CLM-jfa2017}.

\begin{lem}\label{neg0}
Let $\mu$ be a Borel measure in $\rn$ so that $\mu(K_s)$ is second order differentiable at $s=0$
for any $K\in \kpn$ and $\psi\in C_e^2(\sn)$.
Then $\mu$ is log-concave on origin-symmetric convex bodies if and only if
\begin{equation}\label{neg}
\frac{d^2}{ds^2}\log\mu(K_s)\Big|_{s=0} \le 0,
\end{equation}
for any $K\in \kpn$ and $\psi\in C_e^2(\sn)$.
\end{lem}

\begin{proof}
Assume that $\mu$ is log-concave. When $s, t$ are small, we have $K_s, K_t\in \kpn$ and
\[
K_{(1-\lambda) s+ \lambda t} = (1-\lambda) K_s + \lambda K_t.
\]
Then the log-concavity of $\mu$ gives that
\[
\log \mu(K_{(1-\lambda) s+ \lambda t}) \ge (1-\lambda) \log\mu(K_s) + \lambda \log\mu(K_t),
\]
that is, $\log\mu(K_s)$ is concave and \eqref{neg} follows.

Conversely, for $K, L\in \kpn$, let $K_\lambda = (1-\lambda)K + \lambda L$, and $\psi = h_L -h_K$.
Then $\log \mu(K_\lambda)$, $\lambda \in [0,1]$, is concave. This is because of the following,
\[
\frac{d^2}{d\lambda^2}\log\mu(K_\lambda) = \frac{d^2}{ds^2}\log\mu(K_{\lambda+s})\Big|_{s=0}
 = \frac{d^2}{ds^2}\log\mu((K_{\lambda})_s)\Big|_{s=0} \le 0.
\]
Thus, \eqref{l-c} is true.
\end{proof}

In the following, we consider the measure $\mu$ given by $d\mu = e^{-W} dx$, where
$W: \rn\setminus\{0\} \to \R$ is differentiable.

The second order derivative in \eqref{neg} was calculated by Kolesnikov and Milman \cite{KM-ajm2018},
\begin{align}
\frac{d^2}{ds^2}&\log\mu(K_s)\Big|_{s=0} \nonumber \\
&= \frac1{\mu(K)} \int_{\partial K} \big(H_W \varphi^2 -
\langle \text{II}^{-1}\nabla_{\partial K} \varphi, \nabla_{\partial K} \varphi\rangle\big) d\mu_{\partial K}
-\frac1{\mu(K)^2} \Big(\int_{\partial K} \varphi d\mu_{\partial K}\Big)^2, \label{KM1}
\end{align}
where II is the second fundamental form of $\partial K$,
$H_W = \text{tr(II)} - \langle \nabla W, \nu\rangle$ the weighted mean curvature, and
$d\mu_{\partial K} = e^{-W} d\Hn_{\partial K}$.

The following integral formula was proved by Kolesnikov and Milman in \cite{KM-jfa2017},
\begin{align}
\int_K(\Delta u &-\langle \nabla W, \nabla u\rangle)^2\, d\mu =
\int_K \big(\|\nabla^2u\|^2_2 + \langle \nabla^2 W \nabla u, \nabla u\rangle \big)\, d\mu \nonumber \\
&+\int_{\partial K} (H_W u_\nu^2 - 2\langle \nabla_{\partial K} u, \nabla_{\partial K} u_\nu\rangle +
\langle \text{II}\nabla_{\partial K} u, \nabla_{\partial K} u\rangle)\, d\mu_{\partial K}, \label{KM2}
\end{align}
where $u\in C^2(K)$ and $u_\nu = \langle \nabla u, \nu\rangle$.

It was shown by Kolesnikov and Milman \cite{KM-ajm2018} that the Newmann boundary problem,
\begin{equation}\label{KM3}
	\begin{cases}
		\Delta u -\langle \nabla W, \nabla u\rangle =1,& \ \ \text{in } K, \nonumber \\
		\langle \nabla u, \nu\rangle = f,& \ \ \text{on } \partial K.
	\end{cases}
\end{equation}
has a solution if
\[
\int_K f \, d\mu_{\partial K} = \mu(K).
\]

The following lemma is the limiting case of results in Kolesnikov-Livshyts \cite{kolesnikov-livshyts}
and was shown implicitly. For readability, we follow the arguments with slight changes.

\begin{lem}\label{KL}
If for every $K\in \kpn$ and every smooth $u : K \to \R$ with
$\Delta u -\langle \nabla W, \nabla u\rangle =1$ in $K$,
\[
\int_K \big(\|\nabla^2u\|^2_2 + \langle \nabla^2 W \nabla u, \nabla u\rangle \big)\, d\mu \ge 0,
\]
then $\mu$ is log-concave on origin-symmetric convex bodies.
\end{lem}

\begin{proof}
By Lemma \ref{neg0} and \eqref{KM1}, it suffices to show the inequality,
\begin{equation}\label{neg1}
\int_{\partial K} \big(H_W \varphi^2 -
\langle \text{II}^{-1}\nabla_{\partial K} \varphi, \nabla_{\partial K} \varphi\rangle\big) d\mu_{\partial K}
-\frac1{\mu(K)} \Big(\int_{\partial K} \varphi\, d\mu_{\partial K}\Big)^2 \le 0.
\end{equation}

First, assume $\displaystyle \int_K \varphi \, d\mu_{\partial K} \neq 0$. Since the inequality
\eqref{neg1} is homogeneous about $\varphi$, we can replace $\varphi$ by $c\varphi$ so that
$\displaystyle
\int_{\partial K} c\varphi \, d\mu_{\partial K} = \mu(K). $
Then \eqref{neg1} becomes
\begin{equation}\label{neg2}
\int_{\partial K} \big(H_W \varphi^2 -
\langle \text{II}^{-1}\nabla_{\partial K} \varphi, \nabla_{\partial K} \varphi\rangle\big) d\mu_{\partial K}
-\mu(K) \le 0,
\end{equation}
with $\displaystyle \int_{\partial K} \varphi \, d\mu_{\partial K} = \mu(K). $

Let $f=\varphi$ and $u$ a solution to the Newmann boundary problem \eqref{KM3}. Then \eqref{KM2} and the assumption
give that
\begin{align*}
\mu(K) &=\int_K(\Delta u -\langle \nabla W, \nabla u\rangle)^2\, d\mu =
\int_K \big(\|\nabla^2u\|^2_2 + \langle \nabla^2 W \nabla u, \nabla u\rangle \big)\, d\mu \nonumber \\
&+\int_{\partial K} (H_W u_\nu^2 - 2\langle \nabla_{\partial K} u, \nabla_{\partial K} u_\nu\rangle +
\langle \text{II}\nabla_{\partial K} u, \nabla_{\partial K} u\rangle)\, d\mu_{\partial K} \\
&\ge \int_{\partial K} (H_W \varphi^2 - 2\langle \nabla_{\partial K} u, \nabla_{\partial K} \varphi\rangle +
\langle \text{II}\nabla_{\partial K} u, \nabla_{\partial K} u\rangle)\, d\mu_{\partial K}\\
&\ge \int_{\partial K} (H_W \varphi^2 -
\langle \text{II}^{-1}\nabla_{\partial K} \varphi, \nabla_{\partial K} \varphi\rangle)\, d\mu_{\partial K},
\end{align*}
Where the last inequality follows from an inequality for the positive-definite matrix $\Pi$, $\iprod{\Pi x}{x} + \iprod{\Pi^{-1}y}{y}\ge 2\iprod{x}{y}.$ This shows \eqref{neg2}.

When $\displaystyle \int_K \varphi \, d\mu_{\partial K} = 0$, where $\varphi = h_L(\nu) - h_K(\nu)$,
let
\[
\varphi_\lambda = h_{\lambda L}(\nu) - h_K(\nu) = \lambda h_L(\nu) - h_K(\nu),
\]
where $\lambda \neq 1$ and $\lambda$ is close to 1. Then
$\displaystyle \int_K \varphi_\lambda \, d\mu_{\partial K} \neq 0$.
Thus, \eqref{neg1} holds when $\varphi$ is replaced by $\varphi_\lambda$.
Taking limit $\lambda \to 1$ gives \eqref{neg1} for $\varphi$.
\end{proof}

The following result was proved by
Cordero-Erausquin and Rotem \cite{cordero-rotem}.

\begin{lem}\label{thm:BM-measures}
Let $\mu$ be the measure on $\RR^{n}$ with density $e^{-W}$,
where $W(x)=w(|x|)$,
and $w:[0,\infty)\to(-\infty,\infty]$ is increasing
such that $t\mapsto w(e^{t})$ is convex.
If $K\in\kpn$ and smooth even function $u:K\to\RR$ satisfy
$\Delta u-\langle\nabla W,\nabla u\rangle=1$ in $K$, then
\[
\int_K\left(\left\|\nabla^{2}u\right\|_{2}^{2}+\left\langle \nabla^{2}W\nabla u,\nabla u\right\rangle \right)\, d\mu \ge \frac{1}{n}.
\]
\end{lem}

Let us apply Lemmas \ref{KL} and \ref{thm:BM-measures} to prove Theorem \ref{thm:BM}.

\begin{proof}[Proof of Theorem \ref{thm:BM}]
Assume $0<q\le n$.
	Recall the alternative definition of dual quermassintegrals in \eqref{eq:def-norm}.
	If we define a new measure $\mu$ by $d\mu =|x|^{q-n}\, dx$, then (canceling out the constant) \eqref{eq:BM-Wj} becomes
	\[\mu(K+L)^{\frac{1}{q}}\ge \mu(K)^{\frac{1}{q}} + \mu(L)^{\frac{1}{q}}.\]
	This measure nearly satisfies the condition of Theorem \ref{thm:BM-measures} with $w(t)=(n-q)\log t$, however, $w$ attains the value $-\infty$. Hence, we introduce a finite approximation $w_\varepsilon(t)=(n-q)\log (t+\varepsilon)$, and its associated measure $\mu_\varepsilon$ with density $e^{-w_\varepsilon(|x|)}$.
By Lemmas \ref{KL} and \ref{thm:BM-measures}, $\mu_\varepsilon$ is log-concave on origin-symmetric convex bodies,
	\[
	\mu_\varepsilon((1-\lambda)K+\lambda L) \ge\mu_\varepsilon(K)^{(1-\lambda)}\mu_\varepsilon(L)^{\lambda}.
	\]
	Taking $\varepsilon \to 0$, we find
		\[
	\mu((1-\lambda)K+\lambda L) \ge\mu(K)^{(1-\lambda)}\mu(L)^{\lambda}.
	\]
	Note that $\mu(tK)=t^q\mu(K)$, as one may see for example using the fact that $\rho_{tK}(x)=t\rho_K(x)$. We can now use the classical trick to get the desired inequality (see, e.g. \cite[Section 1.2]{shiri-book}).
	Denote
	\[
	K_1=\mu(K)^{-\frac{1}{q}}K, \ \ L_1=\mu(K)^{-\frac{1}{q}}L, \ \ \lambda = \frac{\mu(L)^{\frac{1}{q}}}{\mu(K)^{\frac{1}{q}} + \mu(L)^{\frac{1}{q}}}.
	\]
	Since $\mu(K_1)=\mu(L_1)=1$, the homogeneity and the dimension-free Brunn-Minkowski inequality of $\mu$
give that
\begin{align*}
\frac{\mu(K+L)}{\left(\mu(K)^{\frac{1}{q}} + \mu(L)^{\frac{1}{q}}\right)^q} &=	 \mu\left(\frac{K+L}{\mu(K)^{\frac{1}{q}} + \mu(L)^{\frac{1}{q}}} \right) \\
  &=\mu((1-\lambda) K_1 + \lambda L_1) \\
  &\ge \mu(K_1)^{1-\lambda}\mu(L_1)^\lambda \\
  &= 1
\end{align*}
	which concludes the proof.
\end{proof}

\section{Reverse isoperimetric inequality for dual quermassintegrals} \label{sec:reverse}

 In this section, we shall use yet another equivalent definition of the dual quermassintegral.
 For any $q<n$,
 \begin{equation}\label{eq:def-gaussian}
 		\int_{\RR^n} \rho_K^{q}(x)d\gamma_n(x) =c_{n,q} \wt{V}_q(K)
 \end{equation}
 where $\gamma_n$ is the usual Gaussian probability measure on $\RR^n$ and again $c_{n,j}$ is a constant which does not depend on $K$.  This formulation simply follows from the homogeneous extension of the radial function, and again an integration in polar coordinates:
 \begin{align*}
 	\int_{\RR^n} \rho_K^{q}(x)d\gamma_n(x) &= \frac1{(2\pi)^{n/2}}
 \int_{\RR^n} \rho_K^{q}(x)e^{-\|x\|^2/2}dx  \nonumber \\
 	&=\frac1{(2\pi)^{n/2}}\int_{\Sph^{n-1}}\int_0^\infty r^{-q}\rho_K^{q}(u)e^{-r^2/2}r^{n-1}dr du    \\
 	&= \frac n{(2\pi)^{n/2}}\wt{V}_q(K) \int_0^\infty r^{n-q-1} e^{-r^2/2} dr \nonumber
 \end{align*}
 where the last integral converges when $q < n$.

Define
\[
E_{\gamma_n}(K)= \int_{\rn} \log\rho_K(x)\, d\gamma_n(x).
\]
Then
\[
\lim_{q\to 0}\Big(\int_{\RR^n} \rho_K^{q}(x)d\gamma_n(x)\Big)^\frac1q = e^{E_{\gamma_n}(K)}.
\]
We will also use the identity
\begin{align*}
	&E_{\gamma_n}(K) = c_0 E(K) + c_1,\\
	& c_0= \frac1{(2\pi)^{n/2}}\int_0^\infty e^{-\frac{r^2}2} r^{n-1}dr, \
	c_1=-\frac1{(2\pi)^{n/2}}\int_0^{\infty} (\log r) e^{-\frac{r^2}2} r^{n-1}dr.
\end{align*}

 The proof of the inequality in Theorem \ref{thm:reverse-isop} will use an inequality of Schechtman and Schmuckenschl\"ager (see also \cite{schmuckenschlager1999extremal, barthe1998extremal}).
 We state the inequality together with the equality condition as follows:

 \begin{thm}[\cite{schechtman1995concentration}]\label{thm:S-S}
 	Let  $K \subseteq \RR^n$ be an origin-symmetric convex body.
 If the ellipsoid of maximal volume contained in $K$ coincides with the Euclidean unit ball,
  then, for each $t> 0$,
 	\[\gamma_n(\{x\in\rn: \rho_K^{-1}(x) > t \})\ge \gamma_n(\{x\in\rn :\rho_{B_\infty^n}^{-1}(x) > t\}), \]
 with equality if and only if $K$ is a rotation of $B_\infty^n$.
 \end{thm}

The original statement of this theorem did not include the proof of the equality cases. It can be shown by using Barthe's continuous Brascamp-Lieb inequality \cite{barthe-cont}, as we present below, or directly by using the continuous Ball-Barthe lemma as illustrated in Lutwak-Yang-Zhang \cite{LYZ-JDG2004}.

  We now introduce some additional notations and definitions.  Let $\nu$ be an isotropic Borel measure on $\sn$, that is, it satisfies
 \[
 |x|^2 = \int_{\sn} |x\cdot u|^2\, d\nu(u).
 \]

 A cross measure is an even isotropic discrete measure concentrated uniformly on
 $\pm u_1, \ldots, \pm u_n$, where $u_1, \ldots, u_n$ are orthogonal unit vectors.

 Denote by $C$ the cube $B_\infty^n$, that is,
 $C=\{x=(x_1, \ldots, x_n)\in \rn : |x_i|\le 1, i=1,\ldots, n \}$.

\begin{thm}[Barthe's Continuous Brascamp-Lieb Inequality \cite{barthe-cont}]
 Let $\nu$ be an isotropic Borel measure on $\sn$, $F : \R \to (0, \infty)$ a
 positive continuous function, and $f = \chara{[a,b]} F$, $a<b$. Then
 \[
 \int_{\rn} e^{\int_{\sn} \log f(x\cdot u)\, d\nu(u)} dx  \le \Big(\int_a^b F(s)ds\Big)^n,
 \]
 with equality if and only if $\nu$ is a cross measure.
\end{thm}

 \begin{proof}[Proof of Theorem \ref{thm:S-S}]
  Given an isotropic probability measure $\nu$ define the origin-symmetric convex body $Z_\infty$ by it support function
 \[
 h_{Z_\infty}(x) = \sup \{|x\cdot u| : u\in \text{supp}(\nu)\}.
 \]
 For $t>0$, $t Z_\infty^* =\{x \in \rn :  h_{Z_\infty}(x) \le t\}.$
The indicator function of $t Z_\infty^*$ can be written as
\[
\chara{t Z_\infty^*}(x) = e^{\int_{\text{supp}(\nu)} \log \chara{[-t, t]}(x\cdot u)\, d\nu(u)}.
\]
Together with the isotropy of $\nu$ and the continuous  Brascamp-Lieb inequality, we have
\begin{align*}
\gamma_n(t Z_\infty^*) &=\frac1{(2\pi)^\frac n2}\int_{\rn} \chara{t Z_\infty^*}(x) e^{-\frac{|x|^2}2}dx \\
&=\frac1{(2\pi)^\frac n2}\int_{\rn} e^{\int_{\text{supp}(\nu)}\log\big(\chara{[-t, t]}(x\cdot u)e^{-\frac{|x\cdot u|^2}2}\big)d\nu(u)}dx \\
&\le \frac1{(2\pi)^\frac n2}\Big(\int_{-t}^t e^{-\frac{s^2}2} ds\Big)^n\\
&=\frac1{(2\pi)^\frac n2}\int_{\rn} \chara{tC} (x) e^{-\frac{|x|^2}2}dx \\
&=\gamma_n(tC),
\end{align*}
where the equality in the inequality holds if and only if $\nu$ is a cross measure, that is,
$Z_\infty^*$ is a rotation of $C$.

By John's theorem, there exists an isotropic measure $\nu$ for which
\[
B_2^n \subset K \subset Z_\infty^*.
\]
Therefore,
\[
\gamma_n(tK) \le \gamma_n(tZ_\infty^*) \le \gamma_n(tC),
\]
and hence for $t>0$,
\[
\gamma_n(\{x\in\rn : \rho_K^{-1}(x)\le t\}) \le \gamma_n(\{x\in\rn : \rho_C^{-1}(x)\le t\})
\]
with equality if and only if $K=Z_\infty^*$ is a rotation of $C$.  The desired inequality
follows from taking complements.
\end{proof}

 \begin{proof}[Proof of Theorem \ref{thm:reverse-isop}]
 	 By the results of Theorem \ref{thm:S-S} for all $K$ in John's position and all $s>0$
\begin{equation}\label{s-s-r}
 \gamma_n(\{x\in\rn : \rho_K(x)< s\}) \ge \gamma_n(\{x\in\rn : \rho_C(x) < s\}),
\end{equation}
 with equality if and only if $K$ is a rotation of $C$.
 	Taking complements gives,
\begin{equation}\label{eq:S-S-radials}
 		 	\gamma_n(\{x\in \rn : \rho_K(x) \ge s\}) \le \gamma_n(\{x\in\rn : \rho_C(x) \ge s\}),
\end{equation}
 with equality if and only if $K$ is a rotation of $C$.

 	We may now apply the layer cake decomposition of integrals to the left-hand side of \eqref{eq:def-gaussian}, to get
\[
 \int_{\RR^n} \rho_K^{q}(x)d\gamma_n(x)= \int_0^\infty \gamma_n(\{x\in\rn : \rho_K^{q}(x) \ge t\}) dt .
\]
 	Thus, when $0<q<n$, in order to show \eqref{eq:thm-reverse-isop} we may equivalently prove
\[
\int_0^\infty \gamma_n(\{x\in\rn : \rho_K^{q}(x) \ge t\}) dt
\le \int_0^\infty \gamma_n(\{x\in\rn : \rho_{C}^{q}(x) \ge t) dt,
\]
 	or
\[
\int_0^\infty \gamma_n(\{x\in\rn : \rho_K(x) \ge t^{1/q}\}) dt
\le \int_0^\infty \gamma_n(\{x\in\rn : \rho_{C}(x) \ge t^{1/q}) dt,
\]
 This inequality indeed holds, as one can see by integrating  inequality \eqref{eq:S-S-radials} with $t^{1/q}$ in place of $s$, and the equality holds if and only if $K$ is a rotation of $C$.

 When $q<0$, the inequality \eqref{eq:thm-reverse-isop} is equivalent to
 \[
 \int_0^\infty \gamma_n(\{x\in\rn : \rho_K(x) < t^\frac1q\}) \ge
 \int_0^\infty \gamma_n(\{x\in\rn : \rho_C(x) < t^\frac1q\}).
 \]
This follows from integrating both sides of \eqref{s-s-r}.

When $q=0$, by the layer cake theorem once more and Theorem \ref{thm:S-S}, we have
\begin{align*}
\int_{\rn} \log\rho_K\, d\gamma_n &= \int_{\rn} \big((\log\rho_K)^+ - (\log\rho_K)^-\Big)\, d\gamma_n \\
&=\int_0^\infty \gamma_n(\{x\in\rn : (\log\rho_K)^+\ge t\})\, dt -
\int_0^\infty \gamma_n(\{x\in\rn : (\log\rho_K)^-\ge t\})\, dt \\
&=\int_0^\infty \gamma_n(\{x\in\rn : \rho_K(x) \ge e^t\})\, dt
- \int_0^\infty \gamma_n(\{x\in\rn : \rho_K(x) \le e^{-t})\, dt \\
&\le\int_0^\infty \gamma_n(\{x\in\rn : \rho_C(x) \ge e^t\})\, dt
- \int_0^\infty \gamma_n(\{x\in\rn : \rho_C(x) \le e^{-t})\, dt \\
&=\int_{\rn} \log\rho_C\, d\gamma_n
\end{align*}
with equality if and only if $K$ is a rotation of $C$.
\end{proof}

\section{Problems}

The theorems presented in this work do not complete the analogy between the classical Brunn-Minkowski theory and the dual theory. It is evident that the equality cases in Lutwak's conjecture were not proved.
\begin{prob}
	Does the equality of the inequality \eqref{eq:BM-Wj} holds only if $K$ and $L$ are dilates
	of each other?
\end{prob}
We reiterate that this is known when $q\le 1$, as was shown by \cite{xi-zhang} using a geometric argument.

It is of great interest to extend Lutwak's conjecture to all convex bodies. In the case of the Gaussian measure, which is the canonical example in Kolesnikov and Livshyts's argument, a dimensional Brunn-Mikowski inequality does not hold for non-symmetric sets as was shown by Nayar and Tkocz \cite{Nayar-Tkocz2013note}. However, unlike the Gaussian measure, the measure $\mu$ introduced in the proof of Theorem \ref{thm:BM} is $q$-homogeneous, which might results in an affirmative answer to the following question.
\begin{prob}
	Does the inequality \eqref{eq:BM-Wj} holds for arbitrary convex bodies when they are appropriately
	positioned?
\end{prob}

A similar problem can be posed regarding our second result. In general, reverse isoperimetric inequalities have been proved for all convex bodies, with the simplex as an extremal case \cite{ball1991volume, barthe1998extremal,schmuckenschlager1999extremal}. We believe this should also be the case for the dual quermassintegrals.
\begin{prob}
	Generalize Ball's volume ratio inequality of arbitrary convex bodies
	to dual quermassintegrals.
\end{prob}

Finally, while our proofs of both Theorem \ref{thm:BM} and \ref{thm:reverse-isop} rely on the fact that $q\le n$, this is not necessary for the statement of the problems. Indeed, the dual quermassintegrals are well defined even when $q>n$.
\begin{prob}
	Does the inequality \eqref{eq:BM-Wj} holds when $q>n$?
\end{prob}

\begin{prob}
	Does the inequality \eqref{eq:thm-reverse-isop} holds when $q>n$?
\end{prob}

\bibliographystyle{amsplain}
\addcontentsline{toc}{section}{\refname}\bibliography{intrinsic}

\end{document}